\documentclass[review]{elsarticle}

\usepackage{amsmath}  \usepackage{amsthm}
\usepackage{amssymb}  \usepackage{rotate}
\usepackage{mathrsfs}

\newtheorem{thm}{Theorem}[section]
\newtheorem{cor}[thm]{Corollary}
\newtheorem{lem}[thm]{Lemma}
\newtheorem{ex}[thm]{Example}
\newtheorem{pro}[thm]{Proposition}
\newtheorem{defn}[thm]{Definition}
\newtheorem{rem}[thm]{Remark}

\begin{document}

\begin{frontmatter}

\title{Commutative MTL-rings}


\author[mymainaddress]{Samuel Mouchili}

\author[mymainaddress]{Surdive Atamewoue\corref{mycorrespondingauthor}}
\cortext[mycorrespondingauthor]{Corresponding author}
\ead{surdive@yahoo.fr}
\author[mymainaddress]{Selestin Ndjeya}

\author[mysecondaryaddress]{Olivier Heubo-Kwegna}

\address[mymainaddress]{Department of Mathematics, University of Yaounde 1, Cameroon}
\address[mysecondaryaddress]{Department of Mathematical Sciences, Saginaw Valley State University, 7400 Bay Road, University Center MI 48710}

\begin{abstract}
We introduce in this work, the class of commutative rings whose lattice of ideals forms an MTL-algebra which is not necessary a BL-algebra. The so-called class of rings will be named MTL-rings. 
We prove that a local commutative ring with identity is an MTL-ring if and only if it is an arithmetical ring. It is shown that a noetherian commutative ring $R$ with an identity is an MTL-ring if and only if ideals of the localization $R_M$ at a maximal ideal $M$ are totally ordered by the set inclusion. Remarking that noetherian MTL-rings are again BL-rings, we work outside of the noetherian case by considering non-noetherian valuation domains and non-noetherian Pr\"{u}fer domains. We established that non-noetherian valuation rings are the main examples of MTL-rings which are not BL-rings. This leads us to some constructions of MTL-rings from Pr\"{u}fer domains: the case of holomorphic functions ring through their algebraic properties and the case of semilocal Pr\"{u}fer domains through the theorem of independency of valuations. We end up giving a representation of MTL-rings in terms of subdirectly irreducible product.
\end{abstract}

\begin{keyword}
MTL-algebra\sep Arithmetical ring\sep Valuation ring\sep Pr\"{u}fer domain\sep Holomorphic fonctions ring\sep Subdirect irreducible ring
\MSC[2010] 03G25\sep  13Fxx
\end{keyword}

\end{frontmatter}


\section{Introduction}

\indent It is known that artificial intelligence (AI) and fuzzy logic are strongly related. Actually, the role of AI is to make a computer simulate a human being in handling both certainty and uncertainty  informations; and a basic tool for computers in dealing with those kind of informations is logic. Classical logic deals with certain informations whereas many-valued logic (fuzzy logic) known as nonclassical logic,  deals with uncertain and fuzzy informations. Nowadays, nonclassical logic has become a very useful technique for computer science to deal with such (uncertain and fuzzy) informations. One of the most famous nonclassical logic is the MTL-logic which has been introduced in 2001 by Godo and Esteva \cite{GE}. But three years before, in 1998, H\'{a}jek introduced a very general many-valued logic, called Basic Logic \cite{Hajek}. Since a t-norm has a residuum if and only if the t-norm is left continuous, then the Basic Logic is not the most general t-norm-based logic. Indeed the MTL-logic was proved by Jenei and Montagna \cite{JM} to be the most general t-norm based logic with residua.\\
\indent Many logics have been represented as algebras. For instance boolean algebras are algebras that best fit classical logic. Similarly, in nonclassical logics, several algebras are studied such as MV-algebras for Many-Valued logic. G\"{o}del algebras for G\"{o}del logic. BL-algebras for BL-logic. MTL-algebras for MTL-logic.\\
\indent Since rings are very important algebraic tools in mathematics which have very good and nice properties, many mathematicians have started to relate those algebras with the well known classes of rings. We have for example Jensen for Arithmetical rings \cite{JensenAR}, Belluce, Di Nola and Marchioni for MV-rings and Godel-rings \cite{Luka, Godel}, Heubo-Kwegna, Lele, Ndjeya and Nganou for BL-rings \cite{BL}.\\
\indent Given any commutative ring $R$ in which the following condition holds,
 $$  for~ all~  x\in R~  there ~exists~ an~ element ~ r\in R ~ such ~ that ~ x=rx.~\qquad (*),$$
the semiring of ideals of $R$ under the usual operations, is a residuated lattice denoted $\mathcal{I}d(R)$. (In this work, without any assumption, all the rings used, are commutative rings with (*) condition). For example, in 2009, rings $R$ for which $\mathcal{I}d(R)$ is an MV-algebra known as Lukasiewicz rings were introduced in \cite{Luka}. One year later (in 2010), the ring structure for  G\"{o}del algebras known as  G\"{o}del rings were introduced (\cite{Godel}). Recently in 2018, BL-rings were presented in \cite{BL} as rings whose ideals form a BL-algebra. The natural question that one can ask is the following. Are there commutative rings whose lattice of ideals can be equipped with the MTL-algebra's structure without being a BL-algebra? \\
\indent The current paper tries to answer this question in the following steps. We start by giving some preliminaries on MTL-algebras. The next step is to define the notion of MTL-rings by weaking axioms of BL-rings (the divisibility property on BL-rings will be removed). Remarking that noetherian MTL-rings are again BL-rings, we work outside of the noetherian case by considering non-noetherian valuation domains and non-noetherian Pr\"{u}fer domains. And this yields some examples of MTL-rings which are not BL-rings. We present some other constructions of MTL-rings. Before discussing the representation of MTL-rings in terms of subdirectly irreducible product, we will investigate other properties of MTL-rings together with the connection between MTL-rings, arithmetical rings and their localizations.

\section{Preliminaries on MTL-algebras}
In this section, we will give definitions, examples and some properties of MTL-algebras.
\begin{defn}
(MTL-algebra)\\
An algebraic structure $(Mtl,\wedge,\vee,\odot,\rightarrow,0,1)$ of type $(2,2,2,2,0,0)$ is called an MTL-algebra if the following conditions are satisfied:
\begin{itemize}
\item[(1)] $(Mtl,\wedge,\vee,0,1)$ is a bounded lattice,
\item[(2)] $(Mtl,\odot,1)$ is a commutative monoid,
\item[(3)] $x\odot y\leq z$ if and only if $x\leq y\rightarrow z$, $\forall x, y, z\in Mtl$.
\item[(4)] $(x\rightarrow y)\vee (y\rightarrow x)=1$, $\forall x,y\in Mtl$.
\end{itemize}
\end{defn}

Note that the first three properties are those of a \textit{residuated lattice} and the fourth condition, the last one, is known as a prelinearity condition. \\
To make it short, $Mtl$ will denote the universe of an MTL-algebra $(Mtl,\wedge,\vee,\odot,\rightarrow,0,1)$.\\

\begin{defn}
Let $Mtl$ be an MTL-algebra. Then $Mtl$ is  called:
\begin{itemize}
\item[(1)] a BL-algebra if $x\wedge y=x\odot(x\rightarrow y)$, for all $x,y\in Mtl$ (\textit{divisibility condition}). \label{defBL}
\item[(2)] an MV-algebra if $(x\rightarrow y)\rightarrow y=(y\rightarrow x)\rightarrow x$.
\item[(3)] a G\"{o}del algebra if $x\odot x = x$, for all $x\in Mtl$.
\end{itemize}
\end{defn}

\textbf{Notation}:\\
For all $x,y\in Mtl$,
\begin{center}
$x^n:=x\odot\cdots\odot x$; ~~$\neg x:=x\rightarrow 0$;~~ $x\leftrightarrow y:=(x\rightarrow y)\wedge (y\rightarrow x)$;~~ $x\leq y \Leftrightarrow x\vee y := y \Leftrightarrow x\wedge y:=x$;~~ $\neg\neg x:=x$;~~ $x\oplus y:=\neg(\neg x\odot \neg y).$
\end{center}
Since MTL-algebras are lattices, the notions of homomorphism, quotient, subalgebra, product are define. For more details, one can refer to \cite{Sanka}.

Here below are some examples of MTL-algebras.

\begin{ex}
According to the Definition \ref{defBL}, an MTL-algebra $A$ is a BL-algebra if and only if in $A$, the divisibility condition holds, that is: $x\wedge y=x\odot(x\rightarrow y)$. So BL-algebras are MTL-algebras.  \cite{antoneta}.
\end{ex}

\begin{ex}
 \cite{JunTao}. Let $Mtl=[0,1]$. We define
$$x\wedge y=\min\{x,y\},\, x\vee y=\max\{x,y\},\,
$$
$$
x\odot y = \left\lbrace
\begin{array}{cccc}
0,& x\leq \neg y\\
x\wedge y,& x>\neg y
\end{array}
\right.
; \qquad
x\rightarrow y = \left\lbrace
\begin{array}{cccc}
1,& x\leq y\\
(1-x)\vee y,& x>y.
\end{array}
\right.
$$
The structure $(Mtl,\wedge,\vee,\odot,\rightarrow,0,1)$
is an MTL-algebra. But it is not a BL-algebra. Indeed, for $0\neq y<x<\frac{1}{2}, \, x\rightarrow y=(1-x)\vee y=1-x$ and $x\odot(x\rightarrow y)=x\odot(1-x)=0$ (since $x\leq \neg x$). But $x\wedge y=y\neq 0$.
\end{ex}

\begin{pro}
 \cite{JunTao}. In any MTL-algebra $Mtl$, the following properties hold. For all $x,y,z\in Mtl$,
\begin{itemize}
\item[(1)] $x\leq y$ if and only if $x\rightarrow y=1$,
\item[(2)] $x\odot y\leq x\wedge y$,
\item[(3)] $1\rightarrow x=x$,
\item[(4)] $x\rightarrow(y\wedge z)=(x\rightarrow y)\wedge (y\rightarrow z)$,
\item[(5)] $(x\vee y)\rightarrow  z=(x\rightarrow z)\wedge (y\rightarrow z)$,
\item[(6)] $x\leq y \Rightarrow x\odot z \leq y\odot z$,
\item[(7)] $x\rightarrow y=x\rightarrow(x\wedge y)$,
\item[(8)] $x\rightarrow y = (x\vee y)\rightarrow y$,
\item[(9)] $x\leq y\rightarrow x$.
\end{itemize}
\end{pro}

\section{MTL-rings}
Given any commutative ring $R$, the lattice of ideals $\mathcal{I}d(R)=(Id(R),\wedge,\vee,\odot,\rightarrow,\{0\},R)$, of $R$ is a residuated lattice, where\\
$I\wedge J=I\cap J$, $I\vee J=I+J$, $I\odot J=I\cdot J$, $I\rightarrow J=\{x\in R\, :\, xI\subset J \}$, for all $I,J\in Id(R)$ the set of ideals of $R$.\\
It is recalled that $I^\bullet:=I\rightarrow 0$ is the annihilator of $I$ in $R$.\\ For instance, given a commutative ring $R$, an ideal $I$ of $R$ is called an annihilator ideal if $I=J^\bullet$ for some ideal $J$ of $R$.
The previous notations are adopted throughout this work.
\begin{defn}
(MTL-ring) \\
A commutative ring with the (*) condition is called an MTL-ring if it satisfies the prelinearity condition denoted here by (MTL): $(I\rightarrow J)+(J\rightarrow I)=R$, for all $I,J\in Id(R)$.
\end{defn}

\subsection{\textbf{Some properties of MTL-rings}}
There are some important properties in MTL-rings. We investigate some of them in what follows.

\begin{pro}
The (MTL) condition in the ring $R$ is equivalent to each of the following conditions:\\
\textbf{MTL-1}: $(I\cap J)\rightarrow K=(I\rightarrow K)+(J\rightarrow K)$ for all ideals $I, J, K$ of $R$.\\
\textbf{MTL-2}: $I\rightarrow (J+K)=(I\rightarrow J)+(I\rightarrow K)$ for all ideals $I, J, K$ of $R$.
\end{pro}

\textbf{Proof}:
Suppose that the (MTL) condition holds.\\
Let us show \textbf{MTL-1}. We want to that $(I\cap J)\rightarrow K=(I\rightarrow K)+(J\rightarrow K)$ for all ideals $I, J, K$ of $R$.\\
$(\supseteq):$ Let $x\in (I\rightarrow K)+(J\rightarrow K)$, $x=x_1+x_2$ for some $x_1\in I\rightarrow K$ and $x_2\in J\rightarrow K$. Show that $x\in (I\cap J)\rightarrow K$.\\
Let $y\in (I\cap J)$. Then $xy=(x_1+x_2)y=x_1y+x_2y\in K$. So $x\in (I\cap J)\rightarrow K.$ Therefore $(I\rightarrow K)+(J\rightarrow K)\subseteq (I\cap J)\rightarrow K$.\\
$(\subseteq):$ Let $x\in (I\cap J)\rightarrow K$. Show that $x\in (I\rightarrow K)+(J\rightarrow K)$. We are looking for $x_1\in I\rightarrow K$ and $x_2\in J\rightarrow K$ such that: $x=x_1+x_2$.\\
Since
\begin{equation}\label{eq1}
(I\rightarrow K)+(K\rightarrow I)=R
\end{equation}
\begin{equation}\label{eq2}
(J\rightarrow K)+(K\rightarrow J)=R
\end{equation}
\begin{equation}\label{eq3}
\left((I\cap J)\rightarrow K)+(K\rightarrow (I\cap J)\right)=R
\end{equation}
we have:\\
$(I\rightarrow K)+(J\rightarrow K)+(K\rightarrow I)+(K\rightarrow J)=\left((I\cap J)\rightarrow K)+(K\rightarrow (I\cap J)\right)
$.\\
Set $I\rightarrow K= A_1$, $J\rightarrow K= A_2$, $K\rightarrow I= B_1$, $K\rightarrow J= B_2$, $(I\cap J)\rightarrow K= C$ and $K\rightarrow (I\cap J)= D$.\\
We have that: $A_1+A_2+B_1+B_2=C+D=R$ and $A_1+A_2\subseteq C$ and $D\subseteq B_1+B_2$.\\
In general, for all $x\in R$, there exist $x_C\in C, x_D\in D, x_{A_1+A_2}\in A_1+A_2$ and $x_{B_1+B_2}\in B_1+B_2$ such that $x=x_C+x_D=x_{A_1+A_2}+x_{B_1+B_2}$.\\
Let $x_C\in C$, $x_C=x_{A_1+A_2}+x_{B_1+B_2}$, for certain $x_{A_1+A_2}\in A_1+A_2$ and $x_{B_1+B_2}\in B_1+B_2$.\\
Since $A_1+A_2\subseteq C$, we can take $x_{B_1+B_2}=0$. So $x_C=x_{A_1+A_2}$: $C\subseteq A_1+A_2$. Therefore $C=A_1+A_2$, which means that $(I\cap J)\rightarrow K=(I\rightarrow K)+(J\rightarrow K)$.\\
Similarly, let $x_{B_1+B_2}\in B_1+B_2$. Then $x_{B_1+B_2}=x_C+x_D$, for certain $x_C\in C$ and $x_D\in D$.\\
Since $D\subseteq B_1+B_2$, we can take $x_C=0$. So $x_{B_1+B_2}=x_D$: $B_1+B_2\subseteq D$. Therefore $B_1+B_2=D$, which means that $K\rightarrow(I\cap J)=(K\rightarrow I)+(K\rightarrow J)$.\\
Similar reasonning can be applied to show \textbf{MTL-2} (and to show that $(J+K)\rightarrow I=(J\rightarrow I)+(K\rightarrow I)$).\\
For the converse, suppose \textbf{MTL-1}. Taking $K=I\cap J$, we have the MTL condition since $I\rightarrow(I\cap J)=I\rightarrow J$, $J\rightarrow(I\cap J)=J\rightarrow I$ and $(I\cap J)\rightarrow (I\cap J)=R$.\\
Now suppose \textbf{MTL-2}. Taking $I=J+K$, we have the MTL condition since $(J+K)\rightarrow J=K\rightarrow J$, $(J+K)\rightarrow K=J\rightarrow K$ and $(J+K)\rightarrow (J+K)=R$. $\square$\\
From the previous proposition, we have the following important corollary.

\begin{cor}
In an MTL-ring $R$, for all ideals $I,J,K$ of $R$, one has:\\
(MTL-1'): $(I+J)\rightarrow K=(I\cap J)\rightarrow K$ and
(MTL-2'): $I\rightarrow(J\cap K)=I\rightarrow (J+K)$.
\end{cor}

\begin{ex}
Every BL-ring is an MTL-ring. The class of BL-rings is studied in \cite{BL}.\\
A totally ordered PID (Principal Ideal Domain) is an MTL-ring. This is mainly because, for all ideals $I, J$ in an ordered totally PID $R$, $I=(a)\subseteq (b)=J$; which is equivalent to $I\rightarrow J = R$. Therefore $(I\rightarrow J)+(J\rightarrow I)=R$.\\
It can be easily seen that it is also a BL-ring.
\end{ex}

The first characterization of MTL-rings is the following equivalence.
\begin{thm} \label{thmequiv}
The commutative ring $R$ with the (*) condition is an MTL-ring if and only if  $A(R)$ its lattice of ideals is an MTL-algebra.
\end{thm}
\textbf{Proof}:
Let $R$ be a commutative ring with (*) condition. Suppose that $R$ is an MTL-ring. We have to show that $A(R)=(Id(R),\cap,+,\cdot,\rightarrow,\{0\}, R)$ is an MTL-algebra.\\
(1) For all ideals $I\in Id(R)$, we have $I\cap\{0\}=\{0\}$ and $I+R=R$: $A(R)$ is then a bounded lattice.\\
(2) Let us show that $(Id(R),\cdot,R)$ is a commutative monoid. For all ideals $I,J,K \in R$.\\
(i) $I\cdot (J\cdot K)=(I\cdot J)\cdot K$.\\
(ii) $I\cdot R=I$ because of the (*) condition. So $R$ is then the neutral element for the law $\cdot $  in $A(R)$.\\
(iii)\\
$ \begin{array}{llll}
I\cdot J &=&\{\sum\limits_{i=1}^{n}x_iy_i,\,x_i\in I\text{ and }y_i\in J,\,i=1,\cdots,n\}\\
&=&\{\sum\limits_{i=1}^{n}y_ix_i,\,x_i\in I\text{ and }y_i\in J,\,i=1,\cdots,n\}\\
&=&J\cdot I.
\end{array}
$\\
The second equality holds because the ring $R$ is commutative. It is deduced from (i), (ii) and (iii) that $(Id(R),\cdot,R)$ is a commutative monoid.\\
(3) Let us now show that for all ideals $I,J,K$ of $R$, $I\cdot J\subseteq K \Leftrightarrow I\subseteq J\rightarrow K$. \\
For all ideals $I,J,K$ of $R$,\\
$\begin{array}{llll}
I\cdot J\subseteq K & \Leftrightarrow & xy\in K,\, \forall x\in I,\text{ and } \forall y\in J\\
& \Leftrightarrow & I\subseteq J\rightarrow K.
\end{array}
$ \\
It can be seen from (1), (2) and (3) that $A(R)$ is a residuated lattice.\\
(4) Since $R$ is an MTL-ring, $\forall I,J\in Id(R)$, $(I\rightarrow J)+(J\rightarrow I)=R$.\\
This last verified condition allows to confirm that $A(R)$ is an MTL-algebra.\\
Conversely, if $A(R)=(Id(R),\wedge,\vee,\odot,\rightarrow,\{0\}, R)$ is an MTL-algebra, where:\\ $I\wedge J=I\cap J$, $I\vee J=I+J$, $I\odot J=I\cdot J$, $I\rightarrow J=\{x\in R\, :\, xI\subset J \}$, for all $I,J\in Id(R)$ the set of ideals of $R$, then $\forall I, J\in \mathcal{I}d(R),\, (I\rightarrow J)+(J\rightarrow I)=R$. So $R$ is an MTL-ring.
$\square$\\
The following lemma will be usefull in establishing the equivalent condition for a commutative ring $R$ with (*) condition to be an MTL-ring.


\begin{lem}\label{ideals}
Let $R$ be a ring. The following conditions hold for all ideals $I,J$ and $K$ such that $I\subseteq J$ and $I\subseteq K$.
\begin{itemize}
\item[\textbf{(a)}] $I\subseteq \left( I^\bullet\cdot J \right)^\bullet$, $I\subseteq J\rightarrow I$, $I\subseteq J\rightarrow K$, $I\subseteq K\rightarrow J$;
\item[\textbf{(b)}] $\left( J/I \right)^\bullet=(J\rightarrow I)/I$ ;
\item[\textbf{(c)}] $\left( J/I \right)\rightarrow (K/I)=(J\rightarrow K)/I$.
\end{itemize}
The sympbole " / '' means the quotient.
\end{lem}

\textbf{Proof}:
Let $R$ be a ring and let $I,J$ and $K$ be ideals of $R$ such that $I\subseteq J$ and $I\subseteq K$.\\
\textbf{(a)} \\
$
\begin{array}{llll}
\left( I^\bullet\cdot J \right)^\bullet &=& \{x\in R\, /\, x(I^\bullet\cdot J)=\{0\}\}\\
&=& \{x\in R\, /\, xy=0, \, \forall y\in I^\bullet\cdot J\}\\
&=& \{x\in R\, /\, x\sum\limits_{i=1}^{n}z_it_i=0, \, z_i\in I^\bullet,\, t_i\in J,\, i=1,\cdots,n;\, n\in \mathbb{N}^\bullet\}\\
&=& \{x\in R\, /\, \sum\limits_{i=1}^{n}xz_it_i=0, \, z_i\in I^\bullet,\, t_i\in J,\, i=1,\cdots,n;\, n\in \mathbb{N}^\bullet\}\\
\end{array}
$.\\
Then $\forall x\in I, \, xz=0,\, \forall z\in I^\bullet$. So $I\subseteq \left( I^\bullet\cdot J \right)^\bullet$.\\
Let $y\in J$. Then $xy\in I$ since $x\in I$ and $I$ is an ideal. So $I\subseteq J\rightarrow I$. The same reasonning can be applied to show that $I\subseteq J\rightarrow K$ and $I\subseteq K\rightarrow J$.\\\\
\textbf{(b)} \\
$
\begin{array}{llll}
\left( J/I \right)^\bullet &=& \{x+I\in R/I\, : \, (x+I)(y+I)= 0_{R/I},\, \forall y\in J\}\\
&=& (J\rightarrow I)/I
\end{array}
$.\\\\
\textbf{(c)} \\
$
\begin{array}{llll}
\left( J/I \right)\rightarrow (K/I) &=& \{x+I\in R/I\, : \, (x+I)(J/I)\subseteq K/I\}\\
&=& \{x+I\in R/I\, : \, xJ/I\subseteq K/I\}\\
&=& (J\rightarrow K)/I.~~\square
\end{array}
$\\

Let us consider the following condition called $(MTL)^*$  for a commutative ring $R$.\\
$(MTL)^*$: for all ideals $I,J$ of $R$, if $I\cap J=\{0\}$, then $I^\bullet+J^\bullet=R$.

\begin{rem} (a) In every commutative ring $R$, the following properties are true for all ideals $I$ and $J$ in $R$.
\begin{enumerate}
\item $I\rightarrow J= I\rightarrow (I\cap J)$;
\item $(I+J)\rightarrow K=(I\rightarrow K)\cap (J\rightarrow K)$;
\item $I\rightarrow (J\cap K)=(I\rightarrow J)\cap (I\rightarrow K)$.
\end{enumerate}
(b) Every MTL-ring satisfies the $(MTL)^*$ condition. Actually since $I\rightarrow J=I\rightarrow (I\cap J)$ and $J\rightarrow I=J\rightarrow (I\cap J)$, then $I\cap J=\{0\} \Rightarrow I\rightarrow J=I\rightarrow 0=I^\bullet$ and $J\rightarrow I=J\rightarrow 0=J^\bullet$. Therefore $R=(I\rightarrow J)+(J\rightarrow I)=I^\bullet+J^\bullet$.
\end{rem}

The quotient ring of $R$ gives a characterization of MTL-ring as shown below.
\begin{pro}
A ring $R$ is an MTL-ring if and only if every quotient (by an ideal) of $R$ satisfies the $(MTL)^*$ condition.
\end{pro}

\textbf{Proof}:
Suppose that $R$ is an MTL-ring. Let $I$ be an ideal of $R$. Let $J$ and $K$ be ideals of $R$ such that $I\subseteq J$, $I\subseteq K$ and $(J/I)\cap (K/I)=I$. Then $J\cap K=I$. Now \\
$
\begin{array}{llll}
(J/I)^\bullet+(K/I)^\bullet &=& (J\rightarrow I)/I+(K\rightarrow I)/I, \text{ from (b) in the Lemma \ref{ideals} }\\
&=& (J\rightarrow(J\cap K))/I+(K\rightarrow(J\cap K))/I\\
&=& \left((J\rightarrow K)+(K\rightarrow J)\right)/I\\
&=& R/I.
\end{array}
$\\
Then $R/I$ satisfies the $(MTL)^*$ condition.\\
Conversely, suppose that every factor of $R$ satisfies the $(MTL)^*$ condition. Let $I,J$ be ideals of $R$, then $R/(I\cap J)$ satisfies the $(MTL)^*$ condition. Since $\left( I/(I\cap J) \right)\cap \left( J/(I\cap J) \right)=I\cap J$, one should have $\left( I/(I\cap J) \right)^\bullet+\left( J/(I\cap J) \right)^\bullet=R/(I\cap J)$. That is $(I\rightarrow J)/(I\cap J)+(J\rightarrow I)/(I\cap J)=R/(I\cap J)$. therefore $\left( (I\rightarrow J)+(J\rightarrow I) \right)/(I\cap J)=R/(I\cap J)$. It follows that $(I\rightarrow J)+(J\rightarrow I)=R$. So $R$ is an MTL-ring.
$\square$

\begin{rem}
An MTL-ring which is not a BL-ring cannot be a multiplication ring.\\
The reason is that for any MTL-ring $R$ which is not a BL-ring, if $R$ is a multiplition ring, then $R$ will satisfy the divisibility low: $\forall I, J\in Id(R), \, I\cap J=I\cdot(I\rightarrow J)$. Since $R$ is an MTL-ring, it will also satisfy the prelinearity low. So $R$ is an BL-ring which is a contradiction.
\end{rem}

The following result gives the characterization of local MTL-rings.

\begin{pro}\label{pro1}
Let $R$ be a commutative local ring with  identity. The ring $R$ is an MTL-ring if and only if $Id(R)$ is a chain under the inclusion $\subseteq$.
\end{pro}

\textbf{Proof}:
Let $R$ be a commutative MTL-ring with identity which is local. Let $M$ be its unique maximal ideal. Let $I$ and $J$ be ideals of $R$. Suppose that $I$ and $J$ are not comparable. That is neither $I\nsubseteq J$ nor $J\nsubseteq I$ hold.\\
$I\nsubseteq J\Rightarrow I\rightarrow J\neq R$ and $J\nsubseteq I\Rightarrow J\rightarrow I\neq R$.\\
Since $R$ is an MTL-ring, $(I\rightarrow J)+(J\rightarrow I) = R$.\\
The ideal $I\rightarrow J$ is included in $M$ by Zorn's lemma. This is because each ideal is included in a maximal ideal of $R$. But there is only one maximal ideal which is $M$. So $I\rightarrow J\subseteq M$.\\
Similarly,  $J\rightarrow I\subseteq M$.\\
Therefore $R=(I\rightarrow J)+(J\rightarrow I)\subseteq M+M=M$. This implies that $M=R$, which contradicts the fact that $M$ is a maximal ideal (a maximal ideal is different from the whole ring, by definition).\\
So $I\subseteq J$ or $J\subseteq I$.\\
The converse is obvious.
$\square$

\subsection{\textbf{Further properties of ideals in MTL-rings}}
Knowing the characterization of ideals in a ring allows to better understand the ring.\\
The following proposition gives futher properties of ideals in MTL-rings.
\begin{pro}
Let $R$ be a commutative MTL-ring with identity, $M$ a maximal ideal of $R$. \\
Let $I$ and $J$ be ideals of $R$. The following properties hold:
\begin{itemize}
\item[(a)] $S^{-1}I\subseteq S^{-1}J$ or $S^{-1}J\subseteq S^{-1}I$, where $S=R\backslash M$ is a multiplication set. And $S^{-1}R$ is the localization at $M$, denoted by $R_M$. We also have the notations: $S^{-1}I=IR_M$ and $S^{-1}J=JR_M$.
\item[(b)] For all ideal $K$ of $R$,
\begin{itemize}
\item[(i)] $I^{a}J^{b}\subseteq I^{a+b}+J^{a+b}$, for all natural integers $a$ and $b$. So
\item[(ii)] $(I+J)^{n}=I^{n}+J^{n}$, for all natural integer $n$.
\item[(iii)] $(I+J)(I\cap J)=IJ$.
\item[(iv)] $K(I\cap J)=KI\cap KJ$.
\item[(v)] $K+(I\cap J)=(K+I)\cap (K+I)$ \label{pro-5}
\item[(vi)] $K\cap(I+J)=(K\cap I)+(K\cap J)$.
\item[(vii)] $K\rightarrow(I+J)=(K\rightarrow I)+(K\rightarrow J)$.
\item[(viii)] $(I\cap J)\rightarrow K=(I\rightarrow K)+(J\rightarrow K)$.
\end{itemize}
\item[(c)] $\sqrt{I+J}=\sqrt{I}+\sqrt{J}$.
\item[(d)] There is an endomorphism $\phi$ of determinent 1 such that :
$$\begin{array}{cccc}
\phi:& R\times R& \longrightarrow & R\times R\\
&(I,J)&\longmapsto & (I+J, I\cap J)
\end{array}
$$
\item[(e)] Chinese Remainder Theorem: $R/I\times R/J \cong R/(I\cap J)\times R/(I+J)$.
\item[(f)] More generally, for all $n\geq 1$, there is an endomorphism $\phi$ of determinent 1 such that:
$$\begin{array}{cccc}
\phi:& R^{n}& \longrightarrow & R^{n}\\
&I=(I_j)_{1\leq j \leq n}&\longmapsto & \left(\sigma_j\left(I\right)\right),
\end{array}
$$
where $\sigma_j(I)$ is homogenious elementary symmetric polynomial of degree $j$ in $I_k$:
$$\sigma_j(I)=\sum_{k_1<\cdots<k_j}I_{k_1}\cap\cdots \cap I_{k_j}=\bigcap_{k_0<\cdots<k_{n-1}}I_{k_0}+\cdots +I_{k_{n-j}}.$$
\end{itemize}
\end{pro}

\textbf{Proof}:
(a) is already proved straightforward because any local MTL-ring is an ideal chain ring.\\
(b) \\
(i) Let $a,b\in \mathbb{N}^*$. We want to show that for all ideals $I$, $J$ of $R$, $I^aJ^b\subset I^{a+b}+J^{a+b}$.\\
Let $M$ be a maximal ideal of $R$. Set $S=R\backslash R$ a multiplicatively set of $R$.\\
Consider the localization $R_M$ of $R$ at $M$. Since $R$ is an MTL-ring, then $R_M$ is an ideal-chain ring (Proposition \ref{pro1}).\\
Then $S^{-1}I\subseteq S^{-1}J$ or $S^{-1}J\subseteq S^{-1}I$.\\
$
\begin{array}{llll}
S^{-1}I\subseteq S^{-1}J &\Rightarrow& \left( S^{-1}I\right)^{b}\subseteq \left(S^{-1}J\right)^{b}\\
&\Rightarrow& \left( S^{-1}I\right)^{b}\left( S^{-1}I\right)^{a}\subseteq \left(S^{-1}J\right)^{a+b}\\
&\Rightarrow& S^{-1}(I^bJ^a)\subseteq \left(S^{-1}J\right)^{a+b}\\
&\Rightarrow& S^{-1}(I^bJ^a)\subseteq \left(S^{-1}J\right)^{a+b}+\left(S^{-1}I\right)^{a+b}\\
&\Rightarrow& I^bJ^a\subseteq I^{a+b}+J^{a+b}\qquad \text{ by the local-global principle } \\
\end{array}
$\\\\
$
\begin{array}{llll}
S^{-1}J\subseteq S^{-1}I &\Rightarrow& \left( S^{-1}J\right)^{a}\subseteq \left(S^{-1}I\right)^{a}\\
&\Rightarrow& \left( S^{-1}J\right)^{a}\left( S^{-1}I\right)^{b}\subseteq \left(S^{-1}I\right)^{a+b}\\
&\Rightarrow& S^{-1}(J^aI^b)\subseteq \left(S^{-1}J\right)^{a+b}\\
&\Rightarrow& S^{-1}(I^bJ^a)\subseteq \left(S^{-1}J\right)^{a+b}+\left(S^{-1}I\right)^{a+b}\\
&\Rightarrow& I^bJ^a\subseteq I^{a+b}+J^{a+b}\qquad \text{ by the local-global principle } \\
\end{array}
$\\\\
So $I^aJ^a\subseteq I^{a+b}+J^{a+b}$.\\
(ii) Now let $n\in \mathbb{N}^*$. Show that $(I+J)^n=I^n+J^n$.\\
By induction on $n\in \mathbb{N}^*$.\\
For $n=2$,\\

$
\begin{array}{llll}
(I+J)^2 &=& I^2+J^2+IJ+JI\\
	    &=& I^2+J^2+2IJ,\text{ because $R$ is commutative }\\
	    &=& I^2 + J^2 +IJ, \text{ since } I,J \text{ are ideals }\\
	    &\subseteq& I^2+J^2+I^2+J^2 \text{ because } IJ\subseteq I^2+J^2 \text{ from } (i).\\
	    &=& I^2+J^2.
\end{array}
$\\
Conversely, $I\subseteq I+J$ and $J\subseteq I+J$ imply $I^2\subseteq (I+J)^2$ and $J^2\subseteq (I+J)^2$. Therefore $I^2+J^2\subseteq (I+J)^2$. So $(I+J)^2=I^2+J^2$.\\
Suppose that $(I+J)^{n-1}=I^{n-1}+J^{n-1}$ and show that $(I+J)^n=I^n+J^n$. \\
Basically, $I^n+J^n\subseteq (I+J)^n$.\\
Now
$
\begin{array}{llll}
(I+J)^n &=& (I+J)^{n-1}(I+J)\\
        &=& I^n+J^n+IJ^{n-1}+JI^{n-1}\\
        &\subseteq& I^n+J^n+I^n+J^n+I^n+J^n=I^n+J^n.
\end{array}
$\\
(iii) Show that $(I+J)(I\cap J)= IJ$.\\
$
\begin{array}{llll}
TIJ=[(I\rightarrow J)+(J\rightarrow I)]IJ&=&I(I\rightarrow J)J+J(J\rightarrow I)I, \text{ since $R$ is commutative}\\
        &\subseteq& (I\cap J)J+(J\cap I)I \\
        &=& (I+J)(I\cap J)\\
        &\subseteq& IJ+IJ=IJ, \text{ since $R$ is commutative}
\end{array}
$\\
(iv) Show that $K(I\cap J)=KI\cap KJ$.\\
$
\begin{array}{llll}
T(KI\cap KJ)&=&\left[\left(I\rightarrow J\right)+\left(J\rightarrow I\right)\right](KI\cap KJ)\\
&\subseteq& (KI)(I\rightarrow J) + (KJ)(J \rightarrow I)\\
&\subseteq& K(I\cap J)+K(J\cap I)=K(I\cap J)\subseteq KI\cap KJ.
\end{array}
$\\
(v) Show that $K+(I\cap J)=(K+I)\cap (K+J)$.\\
$
\begin{array}{llll}
T\left[\left(K+I\right)\cap\left(K+J\right)\right]&=& \left[\left(I\rightarrow J\right)+\left(J\rightarrow I\right)\right]\left[\left(K+I\right)\cap\left(K+J\right)\right]\\
&=& (I\rightarrow J)((K+I)\cap(K+J))+(J\rightarrow I)((K+I)\cap(K+J))\\
&\subseteq& (I\rightarrow J)(K+I)+(J\rightarrow I)(K+J)\\
&=& K(I\rightarrow J)+I(I\rightarrow J)+K(J\rightarrow I)+J(J\rightarrow I)\\
&\subseteq& (I\cap J)+K(I\rightarrow J)+(J\cap I)+K(J\rightarrow I)\\
&\subseteq& K[\underbrace{(I\rightarrow J)+(J\rightarrow J)}_{R}]+(I\cap J)\\
&=& K+(I\cap J)\\
&\subseteq& (K+I)\cap(K+J).
\end{array}
$\\
(vi) Show that $K\cap(I+J)=(K\cap I)+(K\cap J)$.\\
$
\begin{array}{llll}
T(K\cap(I+J))&=&[(I\rightarrow J)+(J\rightarrow I)][K\cap(I+J)]\\
&=& (I\rightarrow J)(K\cap(I+J))+ (J\rightarrow I)(K\cap(I+J))\\
&\subseteq& K(I\rightarrow J)\cap (I\rightarrow J)(I+J)+K(J\rightarrow I)\cap (J\rightarrow I)(I+J)\\
&\subseteq& K(I\rightarrow J)\cap [(I\cap J)+J(I\rightarrow J)]+ K(J\rightarrow I)\cap [I(J\rightarrow I)+(I\cap J)]\\
&\subseteq& K\cap J + K\cap I\\
&\subseteq& K\cap(I+J).
\end{array}
$\\
(vii) Show that $K\rightarrow (I+J)=(K\rightarrow I)+(K\rightarrow J)$.\\
$\begin{array}{llll}
T(K\rightarrow (I+J))&=&[(I\rightarrow J)+(J\rightarrow I)][K\rightarrow (I+J)]\\
&=& (I\rightarrow J)(K\rightarrow (I+J))+(J\rightarrow I)(K\rightarrow(I+J))\\
&\subseteq& (K\rightarrow J)+(K\rightarrow I)\\
&\subseteq& K\rightarrow(I+J).
\end{array} $\\
(viii) Show that $(I\cap J)\rightarrow K=(I\rightarrow K)+(J\rightarrow K)$.\\
$
\begin{array}{llll}
((I\cap J)\rightarrow K)T&=&[(I\cap J)\rightarrow K][(I\rightarrow J)+(J\rightarrow I)]\\
&=& ((I\cap J)\rightarrow K)(I\rightarrow J)+((I\cap J)\rightarrow K)(J\rightarrow I)\\
&\subseteq& (I\rightarrow K)+(J\rightarrow K)\\
&\subseteq& (I\cap J)\rightarrow K.
\end{array}
$\\
(c) Show that $\sqrt{I+J}=\sqrt{I}+\sqrt{J}$.\\
$
\begin{array}{llll}
T\sqrt{I+J}&=&[(I\rightarrow J)+(J\rightarrow I)]\sqrt{I+J}\\
&=& (I\rightarrow J)\sqrt{I+J}+(J\rightarrow I)\sqrt{I+J}\\
&\subseteq& \sqrt{I+J}.
\end{array}
$\\
(d) Let $(i,j)\in (I\rightarrow J)\times (J\rightarrow I)$ such that $i+j=1$. One can define $\phi$ to be the left multiplication by the matrix
$\begin{pmatrix}
1&1\\-i&j
\end{pmatrix}
$\\
$
\begin{array}{llll}
\phi:& R\times R &\longrightarrow & R\times R\\
&(a,b)&\longmapsto &\begin{pmatrix}
1&1\\-i&j
\end{pmatrix} \begin{pmatrix}
a\\b
\end{pmatrix}=\begin{pmatrix}
a+b\\-ai+bj
\end{pmatrix}.
\end{array}
$\\
(e) By the Chinese Remainder Theorem. Just notice that from (d) one has $(I+J)(I\cap J)=IJ$.\\
(f) comes from (d). Repeat (d) to construct $\phi$ and use distributivity: $I\cap (J+K)=(I\cap J)+(I\cap K)$.
$\square$

\section{Link between MTL-rings, Arithmetical rings and localization of rings}
We start this section with some recalls on arithmetical rings.

\begin{defn}\cite{JensenAR}. \label{defJen}
An arithmetical ring $R$ is a commutative ring  with identity whose lattice of ideals is distributive. That is, for all ideals $I, J, K$ of $R$,
$(I+J)\cap K=(I\cap J)+(I\cap K)$.\\
For a commutative local ring with identity, Jensen \cite{JensenAR} showed that it is arithmetic if and only if its lattice of ideals forms a chain under the set inclusion.
\end{defn}

\begin{rem}\label{rem}
Theorem 3 in \cite{JensenAR} shows that the following conditions are equivalent with $R$ a commutative ring with identity:\\
(a) $R$ is arithmetical.\\
(b) $(I+J)\rightarrow K=(I\rightarrow K)+(J\rightarrow K)$ for arbitrary ideals $I$ and $J$ and any finitely generated ideal $K$.\\
(c) $I\rightarrow (J\cap K)=(I\rightarrow J)+(I\rightarrow K)$ for finitely generated ideals $J$ and $K$, and arbitrary $I$.
\end{rem}

If the conditions (b) or (c) is satisfied for any ideals taken arbitrary, then arithmetical rings would be MTL-rings. However, from the Remark \ref{rem} and the (MTL-1') condition,  we have the following proposition.

\begin{pro}\label{pro}
An MTL-ring with an identity is an arithmetical ring.
\end{pro}
The converse of the Proposition \ref{pro} is true only for noetherian rings. Indeed if a non-noetherian arithmetical ring is an MTL-ring, then the prelinearity condition will hold for all ideals, which means that all ideals are finitely generated since in an arithmetical ring, that condition of prelinearity hold only for finitely generated ideals, and this yields a contradiction.

\begin{lem}\cite{JensenAR}.
Let $R$ be a commutative ring with identity. Let $S$ be a multiplicatively closed set not containing zero. If $I$ and $J$ are ideals of $R$ such that $J$ is finitely generated, then $(I\rightarrow J)_S=I_S\rightarrow J_S$, where $I_S$ and $J_S$ are ideals in the localization $R_S$.
\end{lem}

\begin{thm}
Let $R$ be a noetherian commutative ring with identity. $R$ is an MTL-ring if and only if ideals of a localization $R_M$ at any maximal ideal $M$ are totally ordered by the set inclusion.
\end{thm}

\textbf{Proof}:
To prove the implication, suppose that $R$ is an MTL-ring. Then $R_M$ the localization of $R$ at the maximal ideal $M$ is a local MTL-ring, and from the Proposition \ref{pro1}, it is a chain, that is, its ideals are totally ordered by the set inclusion.\\
The converse is obvious because any ideal chain ring is an MTL-ring.
$\square$

\begin{pro} \cite{Luka} \label{proLuka}
Let $R$ be a commutative ring with identity such that\\
1) $R$ is local and Artinian with maximal ideal $M$;\\
2) the lattice of ideals of $R$ is a chain under $\subseteq$.\\
Then $R$ is a Lukasiewicz ring.
\end{pro}

From the Proposition \ref{proLuka}, we find out the following corollary.

\begin{cor}
A noetherian local commutative MTL-ring with a unit is a Lukasiewicz ring.
\end{cor}

The following proposition gives the strong relation between MTL-ring and BL-ring.

\begin{pro}\label{pro2}
Let $R$ be a noetherian commutative ring with an identity. The ring $R$ is an MTL-ring if and only if it is a BL-ring.
\end{pro}

\textbf{Proof}:
Let $R$ be a noetherian commutative ring with identity.\\
$(\supseteq)$ if $R$ is a BL-ring, then $R$ is an MTL-ring from the definition of BL-rings.\\
$(\subseteq)$ If $R$ is an MTL-ring, then $R$ will satisfy the conditions (MTL-1) and (MTL-1'). Then $R$ is arithmetical ring by the Theorem 3 in \cite{JensenAR}. Since $R$ is noetherian, then $R$ is a multiplication ring (according to the Theorem 2 in \cite{JensenAR}). So $R$ is a BL-ring.
$\square$

\begin{rem}
According to Proposition \ref{pro2}, if we intend to find an MTL-rings which is not also a BL-ring, we need to avoid noetherian rings.
\end{rem}

\section{Some contructions of MTL-rings which are not BL-rings}

\begin{defn}\cite{ValRing}.
A commutative ring $R$ with identity whose ideals are totally ordered by inclusion is called $valuation\, ring$. A valuation ring without non-zero zero divisors is called a $valuation\, domain$.
\end{defn}

\begin{ex}
Let $\mathbb{Z}$ be usual ring of integers. The localization $\mathbb{Z}_p$ of $\mathbb{Z}$ at a prime ideal $p\mathbb{Z}$ is a valuation domain. It is its $p$-adic completion $J_p$ (the ring of $p$-adic integers). Since its unique maximal ideal is principal, then all its non-zero ideals are principal: it is a PID(Principal Ideal Domain). The quotient ring $J_p/p^nJ_p=\mathbb{Z}_p/p^n\mathbb{Z}_p\cong \mathbb{Z}/p^n\mathbb{Z}$ is a valuation ring with a finite number of ideals.
\end{ex}

Note that a valuation ring $R$ contains a unique maximal ideal that will be denoted by $P$, that is, valuation rings are local rings. The ideal $P$ is the set of all non-invertible elements of $R$. The field $R/P$ is called the residue field of $R$.

A first link between valuation rings and MTL-rings are stated in the following proposition.

\begin{pro}\label{pro3}
Valuation rings are MTL-rings.
\end{pro}
\textbf{Proof}:
Let $A$ and $B$ be ideals of a valuation ring.\\
Then $A\subseteq B$ or $B\subseteq A$. So $(A\rightarrow B)+(B\rightarrow A)=R$.
$\square$

The following remark is very important.

\begin{rem}
Note that, if in the Proposition \ref{pro3}, we assume that those valuation rings are also noetherian, then these MTL-rings will be BL-rings [Proposition \ref{pro2}]. That is the reason why non-noetherian valuation rings will get our attention to find MTL-rings which are not BL-rings.
\end{rem}

Since multiplication rings correspond to the divisibility property in BL-rings, the natural question that one should ask is the following: can a non-noetherian valuation ring be a multiplication ring? We answer it in what follows.

\subsection{\textbf{Non-noetherian valuation rings}}

\begin{defn}\cite{DDA-MRing}.
Let $R$ be a commutative ring with an identity. An ideal $A$ of $R$ is called a \textit{multiplication ideal} if for every ideal $B\subseteq A$ there exists an ideal $C$ such that $B=AC$. A ring $R$ is called a \textit{multiplication ring} if all its ideals are multiplication ideals.
\end{defn}

\begin{rem}\label{rem1}
(Theorem 1, \cite{DDA-MRing})\\ In a quasi-local ring, every multiplication ideal is principal.\\
The Corollary 2.2 in \cite{DDA-MRing} mentioned that a semi-quasi-local multiplication ring is principal.
\end{rem}
This Remark \ref{rem1} yields the following important theorem.

\begin{thm}\label{theo}
A non-noetherian valuation ring is an MTL-ring which is not a BL-ring.
\end{thm}

\textbf{Proof}:
Let $R$ be a non-noetherian valuation ring. The first part of the proposition is done in the Proposition \ref{pro3}. That is, a valuation ring is an MTL-ring.\\
For the second part, we have to show that $R$ is not a BL-ring. It is shown in \cite{BL} that the divisibility property of BL-rings is equivalent to the one which makes a ring a multiplication ring. Now suppose that $R$ is a BL-ring. Then $R$ is a multiplication ring. Since $R$ is local, then $R$ is principal. This contradicts the fact that $R$ is a non-noetherian ring.
$\square$\\

%
%
Recall that every valuation domain $R$ is the valuation ring $R_v$ of a valuation $v$ of its quotient field (Theorem 3.1 in \cite{ValRing}). According to the same authors Fuchs et al., if $\Gamma(R)\cong\mathbb{Z}$, then a valuation domain $R$ is noetherian. In \cite{AG}, Andreas Gathmann, ensures that the converse is also true.

\begin{ex} Some example of valuations rings which are MTL-rings.
\begin{itemize}
\item[(a)] \textbf{Puiseux series}. Set $\mathcal{P}=\bigcup_{n\in\mathbb{N}}K\left(\left(t^{\frac{1}{n}}\right)\right)$. It is the field of \textit{Puiseux series} on $K$. Let us define on $\mathcal{P}$ a valuation which extends the one define on $K((t))$, the field of formal Laurent series. The value group is $\mathbb{Q}$ since $v(t^{\frac{1}{n}})=\frac{1}{n}$. So the valuation ring is not noetherian.
\item[(b)] More generally, it is noticed in  \cite{ValRing} that for a valuation domain with value group $\mathbb{R}$, there are two isomorphy classes of non-zero ideals: one is the class of principal ideals while the other is the class of all non-principal ideals (these are all isomorphic to the not finitely generated maximal ideal $P$ of $R$).
\item[(c)] Take $R=\mathbb{R}[x,y]$, where $\mathbb{R}$ is the field of real numbers. Define the standard valuation $v:\mathbb{R}[x,y]\rightarrow \mathbb{Z}^2$, with $v(x)=(1,0)\leq v(y)= (0,1)$ and take the value of a polynomial as the minimal values among those of its monomials. The value group is $\mathbb{Z}^2$, which has rank 2. The valuation ring $R_v$ is not noetherian. In general, valuation rings of rank $\geq 2$ is not noetherian.
\end{itemize}
\end{ex}

\subsection{\textbf{Construction of MTL-rings from Pr\"{u}fer domains}}


\begin{defn}\cite{ValRing}
An integral domain $R$ is a Pr\"{u}fer domain if all its localizations at maximal ideals are valuation domains.
\end{defn}

\begin{defn}
A ring without zero-divisor is called a B\'{e}zout ring if all ideal of finite type is principal.
\end{defn}

\begin{ex}
1) A valuation ring is a B\'{e}zout ring.\\
2) All B\'{e}zout rings are Pr\"{u}fer rings, \cite{ValRing}.
\end{ex}

%
The following lemma is the theorem of independency of valuations whose proof can be seen in the Appendix (I) of \cite{Nagata} and in \cite{ValRing} (Theorem 1.7). It gives a method for constructing semilocal Pr\"{u}fer domains which are also MTL-rings.

\begin{lem}\label{lemNagata} \cite{ValRing, Nagata}.
Let $V_1,\cdots,V_n$ be valuation domains with the same quotient field $Q$, such that $V_i\nsubseteq V_j$ for $i\neq j$. Let $P_i$ denote the maximal ideal of $V_i\, (i=1,\cdots,n)$. Then the ring $R=\bigcap_{i=1}^nV_i$ is a semilocal Pr\"{u}fer domain with exactly $n$ maximal ideals: $P_1\cap R,\cdots,P_n\cap R$ and $R_{P_i\cap R}=V_i,\, i=1,\cdots, n$.
\end{lem}
We use the Lemma \ref{lemNagata}, to prove the following theorem.
\begin{thm}\label{theoNagata}
Let $V_1,\cdots,V_n$ be non-noetherian valuation domains with the same quotient field $Q$, such that $V_i\nsubseteq V_j$ for $i\neq j$. Let $P_i$ denote the maximal ideal of $V_i\, (i=1,\cdots,n)$. Then the ring $R=\bigcap_{i=1}^nV_i$ is an MTL-ring which is not a BL-ring.
\end{thm}

\textbf{Proof}:
We prove this theorem for $n=2$. Let $R=V_1\cap V_2$, where $V_1$ and $V_2$ are two non-noetherian valuation domains with the same quotient field $Q$, such that $V_1\nsubseteq V_2$. Let $P_1$ be the maximal ideal of $V_1$ and $P_2$ the maximal ideal of $V_2$.\\
Now let $I$ and $J$ be ideals of $R$. We have the following four cases:
\begin{itemize}
\item[($a_1$)] $I,J\subseteq P_1\cap R$, or
\item[($a_2$)] $I,J\subseteq P_2\cap R$, or
\item[($b_1$)] $I\subseteq P_1\cap R$, and $J\subseteq P_2\cap R$ and $J\nsubseteq P_1\cap R$ and $I\nsubseteq P_2\cap R$ or
\item[($b_2$)] $J\subseteq P_1\cap R$, and $I\subseteq P_2\cap R$ and $I\nsubseteq P_1\cap R$ and $J\nsubseteq P_2\cap R$.
\end{itemize}
\underline{\textbf{Case $(a_1)$}}. $I,J\subseteq P_1\cap R\subseteq V_1 \, \Rightarrow I\subseteq J\text{ or } J\subseteq I$ since $V_1$ is a valuation ring. \\Then $(I\rightarrow J)+(J\rightarrow I)=R$.\\\\
\underline{\textbf{Case $(a_2)$}}. This case is similar to the case ($a_1$); just replace $P_1$ by $P_2$ to have the same result.\\\\
\underline{\textbf{Case $(b_1)$}}. $I\subseteq P_1\cap R$, and $J\subseteq P_2\cap R$ and $J\nsubseteq P_1\cap R$ and $I\nsubseteq P_2\cap R$.\\
The subset $(I\rightarrow J)+(J\rightarrow I)$ is an ideal of the commutative ring $R$ with identity. It follows that either
$(I\rightarrow J)+(J\rightarrow I)\subseteq P_1\cap R$ or $(I\rightarrow J)+(J\rightarrow I)\subseteq P_2\cap R$ or $(I\rightarrow J)+(J\rightarrow I)$ is a maximal ideal of $R$ or $(I\rightarrow J)+(J\rightarrow I)=R$.
\begin{itemize}
\item[(i)] If $(I\rightarrow J)+(J\rightarrow I)\subseteq P_1\cap R$, then $J\subseteq J+I\subseteq(I\rightarrow J)+(J\rightarrow I)\subseteq P_1\cap R$. Hence $J\subseteq P_1\cap R$ which is a contradiction because $J\nsubseteq P_1\cap R$.
\item[(ii)] If $(I\rightarrow J)+(J\rightarrow I)\subseteq P_2\cap R$, then $I\subseteq J+I\subseteq(I\rightarrow J)+(J\rightarrow I)\subseteq P_2\cap R$. Hence $I\subseteq P_2\cap R$ which is a contradiction because $I\nsubseteq P_2\cap R$.
\item[(iii)] If the ideal $(I\rightarrow J)+(J\rightarrow I)$ is a maximal ideal of $R$, it is then different from the two maximal ideals of $R$: $P_1\cap R$ and $P_2\cap R$, by making use of the two previous cases: (i) and (ii). Hence $R$ will have exctly three maximal ideals which is another contradiction since $R$ only has exctly two maximal ideals: $P_1\cap R$ and $P_2\cap R$, according to the Lemma \ref{lemNagata}.
\item[(iv)] Finally, it remains $(I\rightarrow J)+(J\rightarrow I)=R$.
\end{itemize}
\underline{\textbf{Case $(b_2)$}}. This case is similar to the case $(b_1)$. We have the same result by interverting the role of $I$ and $J$.\\
So in all the cases, $(I\rightarrow J)+(J\rightarrow I)=R$, that is, $R$ is an MTL-ring.\\
\indent To show that $R$ is not a BL-ring, suppose that $R$ is a BL-ring. Then $R$ is noetherian. Hence the localization $R_{P_1\cap R}$ of $R$ at $P_1\cap R$ is noetherian. But $R_{P_1\cap R}$ is equal to a non-noetherian valuation ring $V_1$, which is a contradiction. So $R$ is not a BL-ring. \\
To complete the proof, we just have to generalize it by induction for all $n\in \mathbb{N}$.
$\square$

\begin{pro}
The ring $R$ constructed in the Theorem \ref{theoNagata} yields an MTL-ring which is not a valuation ring.
\end{pro}

\textbf{Proof}:
If $R$ is a valuation ring, then $R$ is local. But the set $\{P_1\cap R,\cdots,P_n\cap R\}$ of maximal ideals of $R$ is not a singleton, which contradicts the fact that $R$ is a local ring. So $R$ is an MTL-ring without being a valuation ring.
$\square$

%

\begin{lem}\label{lem1}
Let R be a Pr\"{u}fer domain. If $R$ is non-noetherian,  then every localization of $R$ at a maximal ideal not finitely generated is non-noetherian.
\end{lem}

\textbf{Proof}:
Let $\varphi : R \rightarrow R_M$ be a localization morphism, $R$ is a Pr\"{u}fer domain and $M$ is a maximal ideal of $R$ not finitely generated.\\
$R_M = S^{-1}R$, where $S = R\setminus M$.\\
The ideal $S^{-1}M$ is the maximal ideal of $R_M$.\\
By the first isomorphism theorem, the canonical morphism $ \overline{\varphi}: R/ker(\varphi) \rightarrow \varphi(R)$ is an isomorphism.\\
$R/ker(\varphi) = R/O_R$ since $\varphi$ is an injection. So $R/ker(\varphi) \simeq R$.\\
Suppose that $S^{-1}M$ is finitely generated.\\
Since $R$ is a Pr\"{u}fer domain, then $S^{-1}M$ is a principal ideal of $R_M$. It follows that $(S^{-1}M) \cap \varphi(R)$ is a principal ideal of $\varphi(R)$.\\
Set $I_S= (S^{-1}M) \cap \varphi(R)$. The ideal $\overline{\varphi}^{-1}(I_S)$ is a principal ideal of $R/ker(\varphi)$ because $\overline{\varphi}$ is an  isomorphism of rings. Hence by the principle of saturation, $\overline{\varphi}^{-1}(I_S)\sim M$: contradiction since $M$ is not
finitely generated. So $S^{-1}M$ is not finitely generated, which means that the ring $R_M$ is non-noetherian. $\square$

The following theorem can be deduced from the previous lemma.
\begin{thm}\label{theo1}
Any localization at an not finitely generated maximal ideal of a non-noetherian Pr\"{u}fer domain is an MTL-ring which is not a BL-ring.
\end{thm}
\textbf{Proof}:
Suppose that $R_M$, a localization at an not finitely generated maximal ideal $M$ of a non-noetherian Pr\"{u}fer domain $R$ is a BL-ring. Then, according to the Lemma \ref{lem1}, $R_M$ is a non-noetherian valuation domain. From the Theorem \ref{theo}, the ring $R_M$ is not a BL-ring, which is a contradiction.
$\square$

\begin{cor}\label{coro1}
A localization of the ring of entire functions at one of its not finitely generated maximal ideal is an MTL-ring which is not a BL-ring.
\end{cor}

\subsection{\textbf{Algebraic properties of holomorphic function ring}}
\begin{defn}
\cite{RCA}. A function $f: D\rightarrow \mathbb{C}$ defined on an open set $D$ in the complex plane is said to be $holomorphic$ on $D$ if the limit
$$f'(z)=\lim_{h\rightarrow 0}\frac{f(z+h)-f(z)}{h}$$
is defined for all $z\in D$. The value of the limit is denoted by $f'(z)$, or $\frac{df(z)}{dz}$, and is refered to as the derivative of the function $f$ at $z$.
\end{defn}
Let us now recall some algebraic properties of $\mathcal{H}(\mathbb{C})$, which come from Walter Rudin's book \cite{RCA}.

\begin{pro}
Invertible elements in $\mathcal{H}(\mathbb{C})$ are functions with no zero.
\end{pro}


\begin{pro}
$\mathcal{H}(\mathbb{C})$ is an integral domain.
\end{pro}


Note that the field of fractions of $\mathcal{H}(\mathbb{C})$ is called the field of meromorphic functions $\mathcal{M}(\mathbb{C})$. An element $h\in \mathcal{M}(\mathbb{C})$ is written $h=\frac{f}{g}$, where $f\in\mathcal{H}(\mathbb{C})$ and $g\in \mathcal{H}(\mathbb{C})\backslash \{0\}$.

\textbf{Divisibility in $\mathcal{H}(\mathbb{C})$}.\\
Let us equip the set $\mathbb{Z}\cup \{+\infty\}$ firstly with an addition $+$ which extends the usual addition of $\mathbb{Z}$ to $\mathbb{Z}\cup \{0\}$ by the following rule: for all $n\in \mathbb{Z}\cup \{+\infty\}$, $n+(+\infty)=(+\infty)+n=+\infty,$ and secondly with an order $\leq$ defined by $\forall n,m\in \mathbb{Z}\cup \{+\infty\}, n\leq m$ iff $\exists n'\in \mathbb{N}\cup \{+\infty\}$ such that $m=n+n'$. The so-defined order extends the natural order in $\mathbb{Z}$ to $\mathbb{Z}\cup \{+\infty\}$.\\
For all $\mathcal{H}(\mathbb{C})\backslash \{0\}$ a function $\mu_f$ from $\mathbb{C}$ to $\mathbb{N}$ by for all $z\in \mathbb{C}$, $\mu_f(z)$ is the multiplicity order of the zero of $f$ at $z$. The function $\mu_0$ is the constant function whose value is $+\infty$. So for all $f\in \mathcal{H}(\mathbb{C})$, $\mu_f(z)=0$ iff $f(z)\neq 0$ and $\mathcal{Z}(f)=\mu_{f}^{-1}(\mathbb{Z}_{\geq 1})$, where $\mathbb{Z}_{\geq 1}=\{n\in\mathbb{Z}\cup \{+\infty\},\,n\geq 1\}$. Note that $f$ is invertible in $\mathcal{H}(\mathbb{C})$ iff $\mu_f=0$. \\
Then
\begin{equation}\label{eq1}
\forall f,g\in \mathcal{H}(\mathbb{C}),\,\, \mu_{fg}=\mu_f+\mu_g.
\end{equation}
The function $\mu_f$ can be extended to the functions in $\mathcal{M}(\mathbb{C})$: if $f\in \mathcal{M}(\mathbb{C})\backslash\{0\},\, f=\frac{g}{h}$, where $g\in \mathcal{H}(\mathbb{C})$ and $h\in \mathcal{H}(\mathbb{C})\backslash\{0\}$, then $\mu_f=\mu_g-\mu_h$. The relation (\ref{eq1}) shows that the extension is well defined. And from this formalism, one has the two following statements:\\
$\left(\forall f\in \mathcal{M}(\mathbb{C}),\, \mu_f\geq 0 ,\text{ iff }f\in \mathcal{H}(\mathbb{C})\,\right)$; and
$\left(\forall f,g\in \mathcal{H}(\mathbb{C}),\,  f \text{ divides } g \text{ iff } \mu_f\leq \mu_g\right)$.

\textbf{Factoriality of $\mathcal{H}(\mathbb{C})$.}\\
The only irreducible elements in $\mathcal{H}(\mathbb{C})$ are the functions $f_{z_0}:z\mapsto z-z_0,\,z_0\in\mathbb{C}$ and his product with invertible elements in $\mathcal{H}(\mathbb{C})$.\\
Indeed, we need first to show that $f_{z_0}$ is irreducible. If $g$ divides $f_{z_0}$, then $\mu_g\leq \mu_{f_{z_0}}$. Since $\mu_{f_{z_0}}(z_0)=1$ and $\mu_{f_{z_0}}(z)=0$ for $z\neq z_0$, from the previous inequality, we should have $\mu_g=0$ or $\mu_g=\mu_{f_{z_0}}$. If $\mu_g=0$, then g is invertible. If $\mu_g=\mu_{f_{z_0}}$, then $\mu_{\frac{f_{z_0}}{g}}=0$.  So $\frac{f_{z_0}}{g}$ is an invertible element in $\mathcal{H}(\mathbb{C})$.\\
Conversely, let $f$ be an irreducible element in $\mathcal{H}(\mathbb{C})$. Since $f$ is not invertible, there exists $z_0\in \mathbb{C}$ such that $f(z_0)=0$. In particular, from the previous inequalities, $f_{z_0}$ divides $f$. Since $f_{z_0}$ is also irreducible, we should have that $\frac{f}{f_{z_0}}$ is invertible.\\
It is observed that the function $z\mapsto\sin(z)\in \mathcal{H}(\mathbb{C})\backslash\{0\}$ has infinite number of zeros. From that observation, we have the following proposition.

\begin{pro}\label{proFacto}
The ring $\mathcal{H}(\mathbb{C})$ is not factorial.
\end{pro}

That Proposition \ref{proFacto} shows that the ring $\mathcal{H}(\mathbb{C})$ is not a principal ideal domain $(PID)$. Since we are looking for non-noetherian rings in order to construct MTL-rings which are not BL-rings, let us investigate to figure out whether the ring $\mathcal{H}(\mathbb{C})$ is noetherian or not.

\begin{lem}\label{Bezout}
 The ring $\mathcal{H}(\mathbb{C})$ is a B\'{e}zout ring.
\end{lem}

\textbf{Proof}:
For the proof, see (Theorem 15.15 in \cite{RCA}).
$\square$

\begin{pro}
The ring $\mathcal{H}(\mathbb{C})$ is not noetherian.
\end{pro}

\textbf{Proof}:
Suppose that the ring $\mathcal{H}(\mathbb{C})$ is noetherian, that is, all its ideals are finitely generated. Since $\mathcal{H}(\mathbb{C})$ is a B\'{e}zout ring, then all its ideals are principal (generated by only one element). Therefore $\mathcal{H}(\mathbb{C})$ is a principal ideals domain. Hence the ring $\mathcal{H}(\mathbb{C})$ is factorial, which is a contradiction according to the Proposition \ref{proFacto}. So the ring $\mathcal{H}(\mathbb{C})$ is not noetherian.
$\square$


\begin{thm}
Let $\mathcal{H}(\mathbb{C})$ be the ring of holomorphic functions on the complex plane $\mathbb{C}$. If $M$ is a maximal ideal of $\mathcal{H}(\mathbb{C})$ not finitely generated, then $\mathcal{H}(\mathbb{C})_M$, the localization of $\mathcal{H}(\mathbb{C})$ at $M$ is an MTL-ring which is not a BL-ring.
\end{thm}

\textbf{Proof}:
Since $\mathcal{H}(\mathbb{C})$ is a non-noetherian B\'{e}zout domain, it is a Pr\"{u}fer domain. Then from the Theorem \ref{theo1}, the localization $\mathcal{H}(\mathbb{C})_{M}$ of $\mathcal{H}(\mathbb{C})$ at a maximal ideal $M$ not finitely generated is a non-noetherian valuation ring. Therefore, from the Corollary \ref{coro1}, $\mathcal{H}(\mathbb{C})_{M}$ is an MTL-ring which is not a BL-ring.
$\square$\\

In this section, we gave some examples of MTL-rings which are not BL-rings. Since notherian MTL-rings are BL-rings, we had to work outside of the noetherian case.

\section{Representation of MTL-rings}

We focus here on finite direct products, arbitrary direct sums and homomorphic images of MTL-rings, before giving a representation theorem for commutative MTL-rings in terms of subdirectly irreducible product.
\begin{pro}
MTL-rings are closed under:\\
(1) finite direct products,\\
(2) arbitrary direct sums,\\
(3) homomorphic images. \label{homo}
\end{pro}

\textbf{Proof}:
The proof of the two first points (1) and (2) can be seen in \cite{BL}, Proposition 2.12.\\
For (3), let $R$ be an MTL-ring and $I$ be an ideal of $R$. We shall show that $R/I$ is an MTL-ring. Recall that ideals of $R/I$ are of the form $J/I$, where $J$ is an ideal of $R$ containing $I$. Let $J,K$ be ideals of $R$ containing $I$.\\
$
\begin{array}{llll}
\left( \left(J/I\right)\rightarrow \left( K/I\right)\right)+\left( \left(K/I\right)\rightarrow \left(J/I\right) \right) &=& \left(J\rightarrow K\right)/I + \left(K\rightarrow J\right)/I\\
&=& \left((J\rightarrow K)+(K\rightarrow J)\right)/{I}\\
&=& R/I,\, \text{ since } (J\rightarrow K)+(K\rightarrow J)=R.
\end{array}
$
$\square$\\

We recall that an algebra $\mathcal{A}$ is a subdirect product of an indexed family $\left( \mathcal{A}_i \right)_{i\in I}$ of algebras if\\
(i) $\mathcal{A}\leq \prod\limits_{i\in I}\mathcal{A}_i$ and\\
(ii) $\pi_i(\mathcal{A})=\mathcal{A}_i$, for each $i\in I$.\\
An embedding $\alpha: \mathcal{A}\rightarrow \prod\limits_{i\in I}\mathcal{A}_i$ is subdirect if $\alpha(\mathcal{A})$ is a subdirect product of the $\mathcal{A}_i$.\\
An algebra $\mathcal{A}$ is subdirectly irreducible if for every subdirect embedding $\alpha: \mathcal{A}\rightarrow \prod\limits_{i\in I}\mathcal{A}_i$, there is an $i\in I$ such that $\pi_i\circ \alpha: \mathcal{A}\rightarrow \mathcal{A}_i$ is an isomorphism.\\
The famous Birkhoff subdirectly irreducible representation theorem of algebras says that every algebra is isomorphic to a subdirect product of subdirectly irreducible algebras (which are its homomorphic images); \cite{Sanka}.

\begin{rem}
MTL-rings are not in general closed under infinite direct product.
\end{rem}

\begin{ex}
The ring $\mathbb{Z}_4$ is an MTL-ring. But the infinite direct product $R=\prod\limits_{i=1}^{+\infty}\mathbb{Z}_4$ is not an MTL-ring since it does not satisfy the condition $(MTL)^*$.\\
Indeed let us consider the following ideals: $$I_n=\left\langle (\underbrace{1,0,1,\cdots,1}_{\text{n first components}},0,0,0,\cdots) \right\rangle \text{ and } J_n=\left\langle (\underbrace{0,0,0,\cdots,0}_{\text{n first components}},1,0,1,\cdots) \right\rangle.$$
$I_n \cap J_n=\{(0,0,\cdots,0)\}$. But $I_n^*+J_n^*\neq R$. So $R$ is not an MTL-ring.
\end{ex}

\begin{thm}
(A representation theorem for commutative MTL-rings)\\
Every commutative MTL-ring $R$ is a subdirect product of a family $\{R_r:r\in R\backslash \{0\}\}$ of subdirectly irreducible MTL-rings satisfying:\\
1. $\mathcal{A}(R)$ is a subdirect product of $\{\mathcal{A}(R_r):r\in R\backslash \{0\}\}$.\\
2. $\mathcal{A}(R_r)$ is an MTL-algebra with a unique atom.
\end{thm}

\textbf{Proof}:
By the Birkhoff subdirectly irreducible representation theorem of algebras, $R$ is a subdirect product of subdirectly irreducible rings, all of whom are homomorphic images of $R$. In fact, for all $r\in R\backslash \{0\}$, by Zorn's lemma, there is an ideal $I_r$ maximal among those ideals that do not contained $r$. It follows that $\bigcap\limits_{r\neq 0}I_r=\{0\}$, each factor $R/I_r$ is subdirectly irreducible and $R$ is a subdirect product of the family  $\{R/I_r\, :\, r\in R\backslash \{0\}\}$. Moreover, each quotient $R/I_r$ is an MTL-ring. Hence, one can take $R_r$ to be $R/I_r$.\\
(1) To show  that $\mathcal{R}$ is a subdirect product of $\{R/I_r\, :\, r\in R\backslash \{0\}\}$, consider $\alpha:\mathcal{A}(R)\rightarrow \prod\limits_{r\neq 0}\mathcal{A}(R/I_r)$ defined by $\alpha(I)(r)=I+I_r\mod I_r$. From the various definitions, one has: \\
(i) $(I\cap J)+I_r\mod I_r=(I+I_r\mod I_r)\cap(J+I_r\mod I_r)$, for all $r\in R\backslash \{0\}$, for all $I,J\in \mathcal{A}(R)$.\\
(ii) $(I+ J)+I_r\mod I_r=(I+I_r\mod I_r)+(J+I_r\mod I_r)$, for all $r\in R\backslash \{0\}$, for all $I,J\in \mathcal{A}(R)$.\\
(iii) $(I\rightarrow J)+I_r\mod I_r=(I+I_r\mod I_r)\rightarrow(J+I_r\mod I_r)$, for all $r\in R\backslash \{0\}$, for all $I,J\in \mathcal{A}(R)$.\\
(iv) $(I\cdot J)+I_r\mod I_r=(I+I_r\mod I_r)\cdot(J+I_r\mod I_r)$, for all $r\in R\backslash \{0\}$, for all $I,J\in \mathcal{A}(R)$.\\
(v) $I+I_r=R\Rightarrow I=R, \, \forall r\in R\backslash\{0\},\, \forall I\in \mathcal{A}(R)$. Indeed, suppose that $I\neq R$. Then $\exists r\in R\, :\, r\notin I \Rightarrow I\subseteq I_r$ because $I_r$ is a maximal ideal among those ideals which don't contained $r$. So $I+I_r=I_r$. But $I+I_r=R$. It follows that $I_r=R$. Therefore $\{0\}=\bigcap\limits_{r\neq 0}I_r=\bigcap\limits_{r\neq 0}R=R$: contradiction because $R\neq\{0\}$.\\
(2) Since $R_r$ is subdirectly irreducible, then each $\mathcal{A}(R_r)$ is an MTL-algebra with a unique atom.
$\square$

\section{Conclusion}
In this work, we introduced the class of MTL-rings. We started by recalling definitions, usefull properties and examples of MTL-rings. Then important properties have been given.\\
Non-noetherian arithmetical rings are the class of rings which are very closed to MTL-rings. Indeed, MTL-rings are arithmetical rings, the converse is true if the ring is noetherian. MTL-rings and BL-rings are strongly connected. Actually, a noetherian commutative MTL-ring with a unit is a BL-ring.\\
Since the notherian MTL-rings are BL-rings, we had to work outside of the noetherian case in order to have examples of MTL-rings which are not also BL-rings.  We have used the result of Anderson in \cite{DDA-MRing} which states that a semilocal multiplication ring is principal. This information leaded us to rings like non-noetherian semilocal Pr\"{u}fer rings, non-noetherian semilocal B\'{e}zout rings, and non-noetherian valuation rings. Making used of algebraic properties of the ring $\mathcal{H}(\mathbb{C})$ of holomorphic functions on $\mathbb{C}$, particularly the fact that $\mathcal{H}(\mathbb{C})$ is a B\'{e}zout ring, we noticed that its localizations at maximal ideals are MTL-rings which are not BL-rings.

\end{document}